\documentclass[reqno,centertags, 12pt]{amsart}
\usepackage{comment,graphicx,url}
\usepackage{color}

\definecolor{purple}{rgb}{0.65, 0, 0.9}
\definecolor{orange}{rgb}{1,.5,0}
\definecolor{gray}{rgb}{0.7,.7,0.7}

\usepackage{enumitem}
\usepackage[]{hyperref}

\usepackage{graphics}
\usepackage{amssymb,amsmath,amsfonts}
-\marginparwidth \oddsidemargin -9mm \evensidemargin \oddsidemargin
\usepackage{latexsym}
\advance\hoffset by 5mm

\def\@abssec#1{\vspace{.1in}\footnotesize \parindent .2in
{\bf #1. }\ignorespaces}
\newtheorem{theorem}{Theorem}[section]

\newtheorem{lemma}[theorem]{Lemma}

\newtheorem{remark}[theorem]{Remark}

\newcommand{\be}{\mathbf e}

\newcommand{\quotes}[1]{``#1''}

\allowdisplaybreaks \numberwithin{equation}{section}

\renewcommand{\be}{\begin{equation}}
\newcommand{\ee}{\end{equation}}

\begin{document}

\title[Growth of curvature/perimeter in the 2D Boussinesq equations]{Growth of curvature and perimeter of temperature patches in the 2D Boussinesq equations}

\author{Jaemin Park}
\thanks{
Departement Mathematik und Informatik, Universit\"at Basel,  Spiegelgasse 1, 4051 Basel, Switzerland jaemin.park@unibas.ch}

\begin{abstract}
In this paper, we construct an example of temperature patch solutions for the two-dimensional, incompressible Boussinesq system with kinematic viscosity such that both the curvature and perimeter grow to infinity over time. The presented example consists of two disjoint, simply connected patches. The rates of growth for both curvature and perimeter in this example are at least algebraic.
\end{abstract}

\maketitle

\section{Introduction and the main results}
In this paper, we investigate the long-time behavior of the two-dimensional incompressible Boussinesq equations  in the absence of thermal diffusivity:
\begin{equation}\label{Boussinesq}
\begin{aligned}
\rho_t + u\cdot \nabla \rho & = 0,\\
u_t + u\cdot \nabla u &= -\nabla p - \rho e_2+ \nu\Delta u, \quad  t>0,\ x\in \mathbb{R}^2\\
\nabla \cdot u&=0,\\
(\rho(t,x),u(t,x))|_{t=0}&=(\rho_0(x),u_0(x)),
\end{aligned}
\end{equation}
where $e_2=(0,1)^T$ and $\nu>0$ is the kinematic viscosity coefficient.  The system \eqref{Boussinesq} describes the evolution of the temperature distribution $\rho$, of a viscous, heat-conducting fluid moving in an external gravitational force field, assuming that the Boussinesq approximation is valid \cite{blandford2008applications,gill1982atmosphere,majda2003introduction,pedlosky1987geophysical}. The primary goal of this paper is to provide an example of a patch-type temperature distribution whose curvature and perimeter grow as time approaches infinity. 

\subsection{Overview of long-time behavior in the Boussinesq equations}
 Before presenting the precise statement of the main theorem, let us provide a brief review of the relevant literature.
 
 \subsubsection{Global well-poseness} In the presence of kinematic viscosity as in \eqref{Boussinesq}, Hou--Li~\cite{HouLi2005} and Chae~\cite{Chae2006} obtained global-in-time regularity results in $(\rho,u)\in H^{m}\times H^{m-1}$ and $(\rho,u)\in H^m\times H^{m}$, respectively, with $m\ge 3$.  In essence, when the initial data are sufficiently smooth and decay rapidly, these results ensure the existence of a global-in-time strong solution.  For cases with rough initial data,  Abidi--Hmidi~\cite{MR2290277} and Hmidi--Keraani \cite{MR2305876} proved the existence of global weak solutions  in $(\rho,u)\in B^0_{2,1}\times \left(L^2\cap B^{-1}_{\infty,1}\right)$ and $(\rho,u)\in L^2\times \left(L^2\cap H^s\right)$, where $s\in[0,2)$ and $B^s_{p,q}$ represents the Besov spaces. 
 
 \subsubsection{Long-time behavior of classical solutions} In the study of long-time behavior of solutions, stability analysis, by itself, provides qualitative information about their long-term behavior, moreover it also plays a crucial role in proving various solution features, as exploited in \cite{choi2021growth,drivas2023twisting}. Denoting 
   \[
   \rho_s^\alpha:=\alpha y,\quad u_s^\beta:=(\beta y, 0)^T,\quad \alpha,\beta\in \mathbb{R},
   \] 
   one can easily see that any pair $(\rho_s^\alpha, u_s^\beta)$ is a steady solution (time-independent) for the system \eqref{Boussinesq}. In \cite{MR4451473}, the authors established stability under perturbations in a Gevrey class near the Couette flow ($(\rho_s^\alpha,u_s^1)$ with $\alpha\le 0$), when considering the Boussinesq system in the spatial domain $\mathbb{T}\times \mathbb{R}$. Near the hydrostatic equilibria ($\rho_s^\alpha,0)$, under perturbations in a Sobolev space,  Doering--Wu--Zhao--Zheng \cite{MR3815212} established stability (when $\alpha<0$) and instability (when $\alpha>0$), considering the Boussinesq system in a general Lipschitz domain. Also Tao--Wu--Zhao--Zheng \cite{tao2020stability} conducted another stability analysis with relaxed assumptions on the initial data in a spatially periodic domain.
   
  Regarding a long-time behavior of solutions without an assumption on the smallness of the initial data,  several quantitative results are available in the literature. Since the density $\rho$ is transported by an incompressible flow, one cannot expect any growth or decay of $\rVert \rho(t)\rVert_{L^p}$. However, a creation of small scale by the flow might induce the growth of finer norms of $\rho$ (or vorticity $\omega:=\nabla \times u$) over time.  Indeed, considering \eqref{Boussinesq} in a bounded domain, Ju \cite{Ju2017} showed  $\rVert \rho \rVert_{H^1}\lesssim e^{ct^2}$, which was further improved to an exponential bound  $e^{ct}$ in $\mathbb{T}^2$ by Kukavica--Wang \cite{KW2020}. Subsequently, Kukavica--Massatt--Ziane \cite{kukavica2021asymptotic} achieved a slightly better upper bound $\rVert \rho \rVert_{H^2}\lesssim C_\epsilon e^{\epsilon t}$ for any small $\epsilon>0$. In addition to these upper bounds on growth rates, an interesting lower bound was obtained by Brandolese--Schonbek \cite{BS2012} proving that in $\mathbb{R}^2$, the kinetic energy $\rVert u\rVert_{L^2}$ must grow faster than $c(1+t)^{1/4}$ as $t\to\infty$, provided that the initial density $\rho_0$ does not have a zero average.  Recently, Kiselev--Park--Yao \cite{kiselev2022small} showed that for a large class of initial data, the Sobolev norms $\rVert \rho \rVert_{H^m}$, for $m\ge 1$, must grow at least at some algebraic rate in $\mathbb{T}^2$ and $\mathbb{R}^2$. Besides these quantitative analyses of norm growth, one might expect some asymptotic behavior due to the damping effect induced by viscosity. In this direction,  Kukavica--Massatt--Ziane \cite{kukavica2021asymptotic} and Aydin--Kukavica--Ziane \cite{aydin2023asymptotic} showed that  $\rVert \nabla u \rVert_{L^2}$ and $\rVert \nu\Delta u - \mathbb{P}(\rho e_2)\rVert_{L^2}$, where $\mathbb{P}$ denotes the Leray projection, converge to $0$ as $t\to \infty$ in a bounded domain.
   

 \subsubsection{Temperature patch problem}
 An interesting class of solutions to a transport equation,
 \begin{align}\label{transport1}
 \rho_t + u\cdot \nabla \rho = 0,
  \end{align}
 is called \textit{patch solutions}. These are weak solutions composed of characteristic functions. For instance, if $\rho_0=1_D$ for some bounded domain $D$, the solution remains as  a characteristic function $\rho(t)=1_{D_t}$ for some time-dependent domain $D_t$, provided that the velocity field $u$ is suitably regular. In a more general sense, in this paper, we refer to a solution $\rho$ as a patch solution, if it can be expressed as a linear combination of characteristic functions defined on some bounded domains.

    As transport  phenomena are prevalent in fluid dynamics, the long-time behavior of patch solutions has been  a subject of active study in various fluid models, particularly concerning the behavior of the patch boundary.  In certain two-dimensional models, it has been observed that the curvature of the patch boundary may grow to infinity (\cite{kiselev2019global} for the 2D Euler), the perimeter also may grow to infinity (\cite{choi2021growth,elgindi2020singular,drivas2023twisting} for the 2D Euler) as $t\to\infty$, or even a singularities can develop  in a finite time (\cite{kiselev2016finite,gancedo2021local} for the generalized surface quasi-geostrophic equations).

   In the context of the 2D Boussinesq equations \eqref{Boussinesq},   the temperature distribution is also transported by the velocity field,  making it natural to explore patch solution $\rho$. The global existence of patch solutions  have been rigorously proved by Danchin--Zhang \cite{danchin2017global} and Gancedo--Garc{\'\i}a-Ju{\'a}rez \cite{gancedo2017global}.  More precisely, \cite[Theorem 4.1]{gancedo2017global} states that    if $\gamma\in (0,1)$ and  $\rho_0=1_{D_0}$ for a simply connected domain $D_0\in\mathbb{R}^2$ with $\partial D_0\in C^{2,\gamma}$, and $u_0\in C^\infty_c(\mathbb{R}^2)$ is divergence-free, then 
   \begin{itemize}
   \item there is a unique global solution  $(\rho(t),u(t))$ such that $u\in L^1((0,T);W^{2,\infty}(\mathbb{R}^2))$ and 
  $ \rho(t)=1_{D_t}$, where $D_t=X_t(D_0)$ with $X_t$, a flow map associated to the velocity field $u$. More precisely $X_t$ is the unique map determine by 
    \begin{align}\label{flow1}
  \frac{d{X}_t}{dt}(x)= u(t,X_t(x)),\quad X_0(x)=x, \text{ for all $x\in \mathbb{R}^2$.}
  \end{align}
  \item  $\partial D_t\in L_{loc}^\infty(\mathbb{R}^+ ; C^{2+\gamma})$,
  \end{itemize}
 which ensures that the curvature cannot grow to infinity in a finite time.\color{black}

 Now, let us consider an initial data $(\rho_0,u_0)$ such that the initial temperature density consists of multiple patches with smooth boundaries. The  existence of global weak solution is guaranteed. Indeed, for a set $D_0$ as described in \ref{assumption1}, it is well-known that $1_{D_0}\in B^{\alpha}_{2,2}$ for any $\alpha<\frac{1}{2}$ (e.g., \cite[Proposition 3.6]{sickel1999pointwise}), thus $\rho_0 \in B^{\alpha}_{2,2}$. In this case, classical embedding theorems in Besov spaces (e.g., \cite[Subsection 2.7]{triebel2010theory}), yield that $\rho_0\in 
 B^0_{2,1} \cap B^0_{p,\infty}$ for $p\in (2,4)$. Combining this with $u_0\in C_c^\infty(\mathbb{R}^2)$, the well-posedness theorem  \cite[Theorem 1.2]{MR2305876} ensures that there exists a unique weak solution $(\rho,u)$  to \eqref{Boussinesq} in the class,
  \begin{align}\label{global_weak}
  (\rho,u)\in C(\mathbb{R}^+; B^0_{2,1} \cap B^0_{p,\infty})\times C(\mathbb{R}^+; H^2).
 \end{align}
  However, technically speaking, the global existence results for patch solutions in \cite{gancedo2017global}  are not directly applicable to the initial data as above, since $\rho_0$ consists of two disjoint patches instead of a single patch. More precisely, it is not trivial to see whether the temperature distribution $\rho$ in \eqref{global_weak} qualifies as  a patch solution, and if it does, whether  boundary regularity can persist; A rough velocity $u\in C(\mathbb{R}^+; H^2)$ does not guarantee enough regularity of a flow map $x\mapsto X_t(x)$ to ensure any regularity of the boundary of the set $X_t(D_0)$.  Since our primary focus  lies elsewhere and considering that the main ideas from \cite{gancedo2017global} can be readily applied to the proof, we will only state  a theorem concerning the global existence of patch solutions involving multiple patches. The detailed proof is left to the interested reader.

 \begin{theorem}\cite[Theorem 3.1]{gancedo2017global}\label{wellposed}
Let $N\in\mathbb{N}$. For $i=1,\ldots,N$, let us  pick real numbers $a_i\in \mathbb{R}$ and simply connected  bounded domains $D_i\subset \mathbb{R}^2$ such that  $\overline{D}_i$ are disjoint and $\partial D_i \in C^{2+\gamma}$. Let us consider initial data $(\rho_0,u_0)$ such that $\rho_0=\sum_{i=1}^N a_i1_{D_i}$ and $u_0\in H^3(\mathbb{R}^2)$ is  divergence-free\color{black}. Then there exists a unique weak solution $(\rho,u)$ to \eqref{Boussinesq} such that $u\in C(\mathbb{R}^+ ; H^2)\cap L^1_{loc}(\mathbb{R}^+; W^{2,\infty}(\mathbb{R}^2))$. In addition,  for almost every $t\ge 0$
 \[
 \rho(t) =\sum_{i=1}^N a_i 1_{D_{i,t}},\quad D_{i,t} = X_{t}(D_i),
  \]
 where  $X_t$ is the flow map generated by the velocity $u$. 
Lastly,  we have  persistence of the curvature, $ \partial D_{i,t}\in L^\infty_{loc}(\mathbb{R}^+; C^{2+\gamma})$ for $i=1,\ldots,N$.\color{black}
 \end{theorem}

\subsection{Main results}
\begin{figure}
\hspace{0.3cm}\includegraphics[scale=1.3]{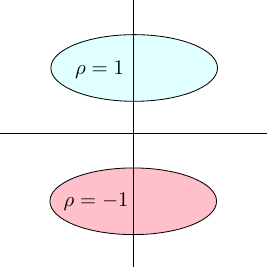}
\caption{Illustration of the initial patch $\rho_0$ }
\label{fig1}
\end{figure}

The goal of this paper is to construct an example of initial data with a patch-type temperature distribution such that under the dynamics \eqref{Boussinesq}, the temperature patch exhibits a growth of the curvature and the perimeter. To this end, we will consider initial data $(\rho_0,u_0)$ satisfying the following assumptions:

\begin{enumerate}[label=\textbf{(A\arabic*)}]
\item\label{assumption1} $\rho_0(x)=1_{D_0}(x) - 1_{D_0^*}$ for a simply connected domain such that $\overline{D_0}\subset \mathbb{R}\times \mathbb{R}^+$, $|D_0|=1$ and $\partial D_0 \in C^\infty$,  and 
\[
D_{0}^*=\left\{ x\in \mathbb{R}^2: (x_1,-x_2)\in D_0\right\}.
\]  
See Figure~\ref{fig1} for an illustration.
\item \label{assumption2} $u_0\in C_c^\infty(\mathbb{R}^2)$ and $u_0$ is divergence-free. Moreover, denoting $u_0=(u_{01},u_{02})$, we assume that  $u_{02}$ is odd in $x_2$ and $u_{01}$ is even in $x_2$.
\end{enumerate}

 When considering $(\rho_0,u_0)$ satisfying \ref{assumption1} and \ref{assumption2}, the uniqueness part of Theorem~\ref{wellposed} ensures  the preservation of   the $x_2$-odd symmetry in the solution $\rho(t)$ obtained in Theorem~\ref{wellposed}: $\rho(t,x_1,-x_2) = -\rho(t,x_1,x_2)$. In other words, $\rho(t)$ takes the form
 \[
 \rho(t) = 1_{D_t} -  1_{D_t^*},\text{ where $D_t=X_t(D_0)$ and $D_{t}^*=\left\{ x\in \mathbb{R}^2: (x_1,-x_2)\in D_t\right\}.$}
 \]
Now, we are ready to state the paper's main theorem:

\begin{theorem}\label{curvature}
Let $(\rho_0,u_0)$ satisfy \ref{assumption1} and \ref{assumption2}. Then, the global patch-type solution $\rho(t)=1_{D_t}-1_{D_t^*}$, where $D_t^*$ denotes the $x_2$ symmetric copy of $D_t$, satisfies the following:
\begin{enumerate}[label=\textbf{(\alph*)}]
\item\label{theorem1} Infinite growth of curvature: We have 
\[
\limsup_{t\to \infty}t^{-\frac{1}{6}}|\kappa(t)| = \infty,
\] 
where $\kappa(t)$ is the maximum curvature of $\partial D_t$.
\item\label{theorem2} Infinite growth of perimeter: Denoting $L_t$ be the distance between a far-left and a far-right point on $\partial D_t$, we have
\[
\limsup_{t\to \infty} t^{-\frac{1}{6}} L_t=\infty.
\]
Since $D_t$ is simply connected, we have infinite growth of perimeter.
\end{enumerate}
\end{theorem}

\section{prelimineries}

  In this section, we collect several well-known conserved properties and some useful uniform estimates for the solutions. 
  
Due to the incompressibility of the flow, the conservation $L^p$ norms of the density follows immediately,
\begin{align}\label{lp_conservation}
\rVert \rho(t)\rVert_{L^p}=\rVert \rho_0\rVert_{L^p},\text{ for all $p\in [1,\infty]$.}
\end{align} 
Another well-known conserved quantity is the total energy of the system. We define the potential energy $E_P(t)$ and the kinetic energy $E_K(t)$ as follows: 
\begin{align*}
E_P(t):=\int_{\mathbb{R}^2}\rho(t,x)x_2 dx,\quad E_K(t):=\frac{1}{2}\int_{\mathbb{R}^2} |u(t,x)|^2dx.
\end{align*}
Then it follows straightforwardly from \eqref{Boussinesq} that
\begin{align}
\frac{d}{dt}E_P(t) &= \int_{\mathbb{R}^2}\rho_t x_2 dx = \int_{\mathbb{R}^2}-u\cdot \nabla \rho x_2 dx = \int_{\mathbb{R}^2}\rho u_2 dx,\label{Epderiv}\\
\frac{d}{dt}E_K(t) &=\int_{\mathbb{R}^2}u\cdot u_t dx = \int_{\mathbb{R}^2}-\rho e_2\cdot u + \nu\Delta u\cdot udx = \int_{\mathbb{R}^2}-\rho u_2 dx - \nu\int_{\mathbb{R}^2}|\nabla u |^2dx.\label{Ekderiv}
\end{align}
By summing up the above two quantities and integrating over time, we obtain
\begin{align}\label{energy_conservation}
E_P(0)+E_K(0) = E_{P}(t)+E_K(t) + \nu\int_0^t \rVert \nabla u(t)\rVert_{L^2}^2dt =:E_T(t)+\nu\int_0^t \rVert \nabla u(t)\rVert_{L^2}^2dt.
\end{align}
Under the assumptions \ref{assumption1} and \ref{assumption2}, both $E_P(t)$ and $E_K(t)$ are always nonnegative. Thus the above energy equality gives us a uniform bound for  the vorticity $\omega:=\nabla\times u$, 
\begin{align}\label{bound2}
\int_0^t \rVert \omega(t)\rVert_{L^2}^2dt \le \int_0^t \rVert \nabla u(t)\rVert_{L^2}^2dt\le C(\rho_0,u_0)(1+{\nu^{-1}}) \text{ for all $t>0$}.
\end{align}

\section{Uniform in time estimates}
In this section, let us derive another uniform time estimate that is simple but will play a crucial role in proving the main theorem. Roughly speaking, the estimate in \eqref{bound2} tells us that time-averaged vorticity dissipates eventually, for instance, $\frac{1}{T}\int_T^{2T}\rVert \omega(t)\rVert_{L^2}^2dt\to 0$, as $T\to \infty$. Now, let us consider the vorticity equation,
\begin{align}\label{vorticity_eq}
\omega_t +u\cdot\nabla \omega = -\partial_1\rho +\nu\Delta\omega,
\end{align}
which can be easily derived by taking the curl operator in the second equation of \eqref{Boussinesq}.  Assuming $\omega$ is sufficiently small in some sense, the quadratic term $u\cdot\nabla \omega$ is comparatively less dominant (in a weak sense) when compared to the linear terms in \eqref{vorticity_eq}. Consequently, considering the dissipation of vorticity, we may anticipate a convergence towards zero for the quantity  $-\partial_1\rho +\nu\Delta \omega$. Establishing such asymptotic behavior in a rigorous sense may  be nontrivial. However, in the next lemma, we will derive a uniform estimate which exhibits  a convergence towards zero of a time-averaged quantity of $-\partial_1\rho +\nu\Delta \omega$.
\begin{lemma}\label{uniform_vs}
Let $(\rho, u)$ be a solution to \eqref{Boussinesq} with initial data $(\rho_0,u_0)$ satisfying  \ref{assumption1} and {\ref{assumption2}}. Then, 
\begin{align}\label{unifo1}
\int_0^t \rVert \partial_1 \Delta^{-1}\rho - \nu\omega\rVert_{\dot{H}^1}^2dt  \le (1+\nu^{-1})C(\rho_0,u_0), \text{ for all $t>0$.}
\end{align}
\end{lemma}
\begin{remark}
When considering \eqref{Boussinesq} in a bounded domain,  in \cite{kukavica2021asymptotic,aydin2023asymptotic}, it was shown that for a general initial data, $\rVert \nu\Delta u - \mathbb{P}(\rho e_2)\rVert_{L^2}$ converges to $0$ as $t\to\infty$, where $\mathbb{P}$ is the Leray projection, which is equivalent to $\rVert\nu \omega -\partial_1\Delta^{-1} \rho\rVert_{\dot{H}^1}\to 0$. When considering the Boussinesq system in an unbounded domain, obtaining such a result for general initial data becomes nontrivial. The challenge arises from the fact that while the proof provided in \cite{kukavica2021asymptotic,aydin2023asymptotic} relies on a uniform bound of the kinetic energy, the kinetic energy in an unbounded domain may not be uniformly bounded in time in general. Although the total energy  is inferred to be bounded from \eqref{energy_conservation}, it remains possible for the kinetic energy to increase indefinitely throughout the evolution without a  lower bound of the potential energy. In our case, this issue is overcome by the assumptions \ref{assumption1} and \ref{assumption2}, which ensure a uniform lower bound of the potential energy.
\end{remark}

\begin{proof}
The proof relies on the second derivative of the potential energy. To begin, we recall the expression of $E''_P(t)$ from \cite{kiselev2022small}:
\begin{lemma}\cite[Lemma 2.1, Lemma 2.3]{kiselev2022small}\label{second_derivative}
Let $(\rho, u)$ be a solution to \eqref{Boussinesq} with initial data $(\rho_0,u_0)$ satisfying \ref{assumption1} and \ref{assumption2}. Then the potential energy $E_P(t)$  satisfies
\begin{equation}\label{2nd_der}
 E_P''(t) = A(t) +   B(t) -\int |\nabla \partial_1\Delta^{-1}\rho(t)|^2dx \quad\text{ for all }t\geq 0,
\end{equation}
where
\begin{equation}\label{def_A}
A(t):=\sum_{i,j=1}^2\int_{\Omega}((-\Delta)^{-1}\partial_2\rho) \partial_i u_j\partial_{j}u_i  \,dx,~~\text{and}~~ B(t):=  \nu  \int_\Omega \rho  \Delta u_2 dx.
\end{equation}
Furthermore, $A(t)$ satisfies
\begin{align}\label{Atbound}
\int_0^t |A(t)|dt \le C(\rho_0)\int_0^t \rVert \nabla u(t)\rVert_{L^2}^2 dt,\text{ for all $t>0$.}
\end{align}
\end{lemma}

Towards the proof of \eqref{unifo1}, we  rewrite $B(t)$, using the Biot-Savart law ($u_2=\partial_1\Delta^{-1}\omega$), as
\[
B(t)=\int \rho \Delta (\partial_1 \Delta^{-1}\omega)dx = -\int \partial_1 \rho \omega dx.
\] Therefore \eqref{2nd_der} can be  also rewritten as
\begin{align}\label{eppdef}
E_P''(t) = A(t) - \nu\int \partial_1 \rho \omega dx - \int |\nabla \partial_1\Delta^{-1}\rho|^2dx .
\end{align}

From the vorticity equation \eqref{vorticity_eq}, we compute
\begin{align}\label{vorticity_energy}
\nu\frac{d}{dt}\left(\frac{1}{2}\rVert \omega(t)\rVert_{L^2}^2\right) &= \nu\int \omega \omega_t dx = \nu \int \omega (-\partial_1 \rho + \nu\Delta \omega)dx = -\int \nu \omega\partial_1\rho dx - \int \nu^2 |\nabla \omega|^2dx.
\end{align}
Combining this with \eqref{eppdef}, we obtain
\begin{align*}
\frac{d}{dt}\left(E'_P(t) + \frac{\nu}2 \rVert \omega(t)\rVert_{L^2}^2\right) &= A(t) -  2\int \nu \omega \partial_1\rho dx - \int |\nabla \partial_1\Delta^{-1}\rho|^2dx  - \int \nu^2 |\nabla \omega|^2dx\\
& = A(t) - 2\int \nabla\nu\omega\cdot \nabla \partial_1\Delta^{-1}\rho dx -  \int |\nabla \partial_1\Delta^{-1}\rho|^2dx  - \int \nu^2 |\nabla \omega|^2dx\\
& = A(t) -  \int |\nabla(\partial_1\Delta^{-1}\rho - \nu \omega)|^2dx.
\end{align*}
Thus integrating this over time, we obtain
\begin{align}\label{uniform_estimates}
E'_P(t) + \frac{\nu}2 \rVert \omega(t)\rVert_{L^2}^2 + \int_0^t \rVert \partial_1 \Delta^{-1}\rho - \nu\omega\rVert_{\dot{H}^1}^2dt =  E_P'(0)+\frac{\nu}2\rVert \omega_0\rVert_{L^2}^2 + \int_{0}^t A(t)dt,\text{ for all $t>0$.}
\end{align}
Finally, sending $E'_P(t)$ on the left-hand side to the other side, 
\[
\int_0^t \rVert \partial_1 \Delta^{-1}\rho - \nu\omega\rVert_{\dot{H}^1}^2dt \le C(\rho_0,u_0) + \int_0^t A(t)dt - E_P'(t) \le (1+\nu^{-1})C(\rho_0,u_0) - E_P'(t),
\]
where the last inequality follows from \eqref{bound2} and \eqref{Atbound}. Recalling $E_P'(t)=\int_{\mathbb{R}^2}\rho u_2dx$ from \eqref{Epderiv}, and using the Cauchy-Schwarz inequality, we have $|E_P'(t)|\le C\rVert \rho(t)\rVert_{L^2}\rVert u\rVert_{L^2}\le C(\rho_0,u_0)$. This gives the desired estimate \eqref{unifo1}.
\end{proof}

\begin{remark}
In the case $\Omega=\mathbb{T}^2$ (or in a suitable bounded domains),  the Sobolev inequality allows us to derive a more concise estimate:
\begin{align}\label{slightlyneater}
\int_0^t\rVert \partial_1 \rho\rVert_{\dot{H}^{-2}}^2dt \le C(\rho_0,u_0)(1+\nu^{-1}).
\end{align}
Indeed, the triangular inequality and the Sobolev inequality give
\[
\rVert \partial_1\Delta^{-1}\rho \rVert_{L^2(\mathbb{T}^2)}^2 \le \rVert \nu\omega \rVert_{L^2(\mathbb{T}^2)}^2 + \rVert \nu\omega - \partial_1\Delta^{-1}\rho\rVert_{L^2(\mathbb{T}^2)}^2 \le \rVert \nu\omega \rVert_{L^2(\mathbb{T}^2)}^2 + \rVert \nu\omega - \partial_1\Delta^{-1}\rho\rVert_{\dot{H}^1(\mathbb{T}^2)}^2.
\]
Thus, integrating over time and combining it with \eqref{unifo1} and \eqref{bound2}, we obtain \eqref{slightlyneater}.
\end{remark}

\section{Lemmas for  curvature and perimeter}
In this section, we study relations between the curvature/perimeter of a patch and the (negative) Sobolev norms. Throughout the section,  $D$ is always assumed to be a simply connected bounded domain such that $\overline{D}\subset \mathbb{R}\times \mathbb{R}^+$, $|D|=1$ and  $\partial D\in C^\infty$. 
 Also,  we will denote the disk centered at $x\in \mathbb{R}^2$ by $B_r(x)\subset \mathbb{R}^2$. A constant $C>0$ will denote a universal constant that does not depend on any variables, while it  might vary from line to line.

Let us recall the Pestov--Ionin theorem, which asserts that every simple-closed curve with a curvature of at most one encloses a unit disk. In other words, it holds that
\begin{align}\label{pIlem}
\sup\left\{r: B_r(x)\subset D,\text{ for some $x\in \mathbb{R}^2$}\right\}\ge  \frac{1}{\max_{x\in\partial D}|\kappa(x)|},
\end{align}
where $\kappa(x)$ is the signed curvature at $x\in \partial D$. In the next lemma, we explore how the radius of a maximal disk within $D$ can be constrained in terms of the negative Sobolev norms of $1_D$.

\begin{lemma}\label{curvature_lem}
Suppose $D$ contains a disk with radius $r>0$.  Then there exists a  universal constant $C>0$ such that 
\begin{align}\label{lemma_cur}
r^3\le C\left( \rVert \partial_1\Delta^{-1}(1_D-1_{D^*}) - \Omega \rVert_{\dot{H}^1(\mathbb{R}^2)} +  \rVert \Omega\rVert_{L^2(\mathbb{R}^2)}\right),
\end{align}
for any function $\Omega$ such that
\begin{align}\label{omega_relaxation}
\rVert \Omega\rVert_{L^2} + \rVert \partial_1\Delta^{-1}(1_D-1_{D^*}) - \Omega \rVert_{\dot{H}^1(\mathbb{R}^2)} <\infty,
\end{align}\color{black}
where $D^*:=\left\{ (x_1,x_2)\subset \mathbb{R}\times \mathbb{R}^-: (x_1,-x_2)\in D\right\}$. 
\end{lemma}  
 The key implication of the above lemma is that the estimate \eqref{lemma_cur} provides a way to measure the maximal curvature of the patch $D$ in terms of a Sobolev norm of $\partial_1\Delta^{-1}(1_D-1_{D^*})$. Intuitively, if $\rVert \partial_1(1_D-1_{D^*})\rVert_{\dot{H}^{-1}}$ is small, then $1_D-1_{D^*}$ must be almost \quotes{independent} of the horizontal variable $x_1$; in other words, the patch $D$ must be horizontally stretched. In such a case, $D$ cannot contain a  large disk in its interior, which, in turn, yields an estimate for the curvature due to \eqref{pIlem}. Our proof of the lemma also lies in the same intuition, combined with a duality of Sobolev norms. A similar idea will be used in the perimeter estimate in Lemma~\ref{periLem}.\color{black}
 \begin{remark}\label{remark_zeromean}
Since $1_D-1_{D^*}\in L^2(\mathbb{R}^2)$ and it has a zero average, we have $\partial_1\Delta^{-1}(1_D-1_{D^*})\in L^2(\mathbb{R}^2)$ \cite[Proposition3.3]{Majda-Bertozzi:vorticity-incompressible-flow}. Therefore, if we  simply take $\Omega:=\partial_1\Delta^{-1}(1_D-1_{D^*})$, then  Lemma~\ref{curvature_lem} tells us that the radius $r$ of a maximal disk contained in $D$ must satisfy
\begin{align*}
r^3\le C\rVert \partial_1\Delta^{-1}(1_D-1_{D^*})\rVert_{L^2(\mathbb{R}^2)}.
\end{align*}
In our proof of the main theorem, we do not have any smallness of $\rVert\partial_1 \Delta^{-1}\rho(t)\rVert_{L^2(\mathbb{R}^2)}$, therefore, we will make a use of carefully chosen $\Omega$ in the application of the lemma.
\end{remark}

\begin{proof}
 Let $B_{0}$ be a disk contained in $D$. Without loss of generality, we  assume that the center of $B_0$ lies on the $x_2$-axis, that is, $B_0=B_r((0,b))$ for some $b>0$. Clearly, we must have $r<1$, due to the assumption $|D|=1$. Next, let us consider a sequence of horizontally translated disks $B_n:=\left\{ (x_1,x_2)\in \mathbb{R}^2: (x_1-2rn,x_2)\in B_r(b)\right\}$ for $n\in\mathbb{N}$ and choose $
 N^*:=\inf\left\{ n \in \mathbb{N}:  |D\cap B_n|\le \frac{r^2}{16}\right\}.$
 We claim that
 \begin{align}\label{nstar}
 |D\cap B_{N^*}|\le \frac{r^2}{16}\text{ and }N^*\le \frac{32}{r^2}.
 \end{align} The first statement is clear by the definition of $N^*$. To see the upper bound of $N^*$, let us suppose, to the contrary, that $N^* > \frac{32}{r^2}$ and denote  $n^*$:=$\lfloor \frac{32}{r^2} \rfloor$, where $\lfloor a\rfloor$ denotes the largest integer not exceeding $a$. Since $r<1$, we have $n^*> \frac{16}{r^2}$. This implies that $|D\cap B_n|\ge \frac{r^2}{16}$ for all $n=1,... n^*$. However, in this case, we must have\[
 1=|D| \ge \sum_{n=1}^{n^*}|D\cap B_n| \ge n^*\frac{r^2}{16}>1,
 \]
 which is a contradiction. 

 Towards the proof of the lemma, we define a function $x_1\mapsto g(x_1)$ and $x_2\mapsto h(x_2)$ as
 \begin{align*}
 g(x_1):=
 \begin{cases}
 0, & \text{ if $x_1\le 0$},\\
 \frac{1-\cos(\pi x_1/r)}{2} & \text{ if $x_1\in (0,r]$},\\
 1 & \text{ if $x_1\in (r,2rN^*]$},\\
 \frac{1+\cos(\pi (x_1-2rN^*)/r)}{2} & \text{ if $x_1\in (2rN^*,2rN^*+r]$},\\
 0 & \text{ if $x_1 > 2rN^*+r$},
  \end{cases}
  \end{align*}
 and
 \begin{align*} h(x_2):= 
 \begin{cases}
 \frac{1+\cos(\frac{\pi (x_2-b)}r)}{2} & \text{ if $x_2\in (b-r,b+r)$},\\
 0 & \text{ otherwise}.
 \end{cases}
 \end{align*}
And we define $f=f(x_1,x_2)$ as 
 \begin{align*}
 \begin{cases}
 f(x_1,x_2):=g(x_1)h(x_2), & \text{ if $x_2\ge 0$,}\\
 f(x_1,x_2):=f(x_1,-x_2),& \text{ if $x_2<0$}.
 \end{cases}
 \end{align*} From the properties of $g$ and $h$, it is clear that the support of $f$ is contained in the vertical strip bounded by $\left\{x_1=0\right\}\cup\left\{ x_1=2rN^*+r\right\}$ whose width is $2rN^* + r\le \frac{C}r$ (see \eqref{nstar}). At the same time, the support of $f$ in $\mathbb{R}\times \mathbb{R}^+$ lies  in the horizontal strip whose width is less than $2r$. Consequently, we have
 \begin{align}\label{support_f}
 |\text{supp}(f)|, \ |\text{supp}(\nabla f)|, \ |\text{supp}(\Delta f)|\le C.
 \end{align}

 Now, denoting $\mu=1_D-1_{D^*}$, we see that for any $\Omega$ satisfying \eqref{omega_relaxation}\color{black},
 \begin{align}\label{l2dual}
 {\int_{\mathbb{R}^2}\mu(x)\partial_1 f(x)dx}&={\int_{\mathbb{R}^2}\mu(x)\Delta^{-1}\Delta\partial_1 f(x)dx}=-{\int_{\mathbb{R}^2}\partial_1\Delta^{-1}\mu(x)\Delta f(x)dx}\nonumber\\
 &= \int_{\mathbb{R}^2} \left(\partial_1\Delta^{-1}\mu(x)-\Omega(x)\right)\Delta f(x)dx +\int_{\mathbb{R}^2} \Omega(x)\Delta f(x)dx\nonumber\\
 &\le \rVert \partial_1\Delta^{-1}\mu - \Omega\rVert_{\dot{H}^1}\rVert \nabla f\rVert_{L^2} +\rVert \Omega\rVert_{L^2}\rVert \Delta f\rVert_{L^2},
 \end{align}
 where we used the  integration by parts and the Cauchy-Schwarz inequality to get the last inequality.
 
 We will estimate the left/right-hand side of the inequality \eqref{l2dual} separately. To get a lower bound of the left-hand side,  we notice that $\mu$ and $\partial_1 f$ are both odd in $x_2$, therefore, 
 \begin{align}\label{numerator1}
 \int_{\mathbb{R}^2}\mu(x)\partial_1 f(x)dx &= 2\int_{\mathbb{R}\times \mathbb{R}^+}\mu(x)\partial_{1}f(x)dx = 2\int_{D}\partial_1 f(x)dx \nonumber\\
 &=2 \int_{D\cap B_0}\partial_1f(x)dx +2 \int_{D\cap B_{N^*}}\partial_1f(x)dx,
 \end{align}
where the last equality follows from the fact that, by the definition of $f$ (especially the definition of $g$),  $\partial_1f(x_1,x_2)=g'(x_1)h(x_2)=0$ if $x\in (B_0\cup B_{N^*})^c$. Using  $B_0\subset D$ and the definition of $g,h$, we can estimate the first integral as
\begin{align*}
\int_{D\cap B_0}\partial_1 f(x)dx &= \int_{B_0} g'(x_1)h(x_2)dx \ge \int_{B_0,\ \left\{\frac{r}{4}<x_1<\frac{3r}{4}\right\}} g'(x_1)h(x_2)dx\\
& \ge \int_{\frac{r}4}^{\frac{3r}{4}}\int_{b-\frac{\sqrt{7}}4r}^{b+\frac{\sqrt{7}}4r} g'(x_1)h(x_2)dx_2 dx_1\\
& \ge \int_{\frac{r}4}^{\frac{3r}{4}}\int_{-\frac{\sqrt{7}}4r}^{\frac{\sqrt{7}}4r} \frac{\pi}{2r}\sin\left(\frac{\pi x_1}r\right) \frac{(1+\cos(\pi x_2/r))}2 dx_2dx_1\\
&=\pi r\int_{1/4}^{1/2}\int_{0}^{\frac{\sqrt{7}}4}\sin (\pi x_1) (1+\cos(\pi x_2)) dx_2dx_1\\
& = \pi r \int_{1/4}^{1/2}\sin(\pi x_1)dx_1\int_0^{\sqrt{7/}4}(1+\cos(\pi x_2))dx_2\\
&\ge \pi r \frac{\sqrt{2}}{2}\frac{\sqrt{7}}4 \ge \pi r\frac{\sqrt{14}}{8}.
\end{align*}
By using $|\partial_1 f|_{L^\infty}\le \frac{\pi}{2r}$, the second integral can be estimated as
\[
\int_{D\cap B_{N^*}}\partial_1f(x)dx\le |D\cap B_{N^*}|\frac{\pi}{2r}\le r^2/16\cdot \frac{\pi}{2r} =\frac{\pi r}{32},
\]
where the second inequality follows from \eqref{nstar}. Thus, in \eqref{numerator1}, we see that 
\begin{align}\label{numerator_estimate}
 \int_{\mathbb{R}^2}\mu(x)\partial_1 f(x)dx \ge 2\left(\frac{\pi r\sqrt{14}}{8} - \frac{\pi r}{32}\right)\ge Cr.
\end{align}

Let us estimate the right-hand side of \eqref{l2dual}   Again, using the properties of $g,h$, we have that $\rVert\Delta f\rVert_{L^\infty}\le \rVert g''\rVert_{L^\infty} + \rVert h''\rVert_{L^\infty}\le \frac{C}{r^2}$ and $\rVert \nabla f \rVert_{L^\infty}\le \rVert g'\rVert_{L^\infty} + \rVert h'\rVert_{L^\infty}\le \frac{C}{r}$. Combining this with \eqref{support_f},  we get 
\begin{align}\label{delta_f}
\rVert \Delta f \rVert_{L^2}\le \frac{C}{r^2},\text{ and }\rVert \nabla f\rVert_{L^2}\le \frac{C}{r}\le \frac{C}{r^2},
\end{align} where the last inequality follows from $r<1$. Plugging this and \eqref{numerator_estimate} into \eqref{l2dual}, we obtain the desired estimate \eqref{lemma_cur}.
\end{proof}

Now, we make a lemma to estimate the parameter of the domain $D$.

\begin{lemma}\label{periLem}
 Let $L>0$ be the distance between a far-left and  a far-right points on $\partial D$, and let $A:=\int_{\mathbb{R}^2}1_{D}(x)x_2dx$.  Then, there exists a  universal constant $C>0$ such that for any $\Omega$ satisfying \eqref{omega_relaxation}\color{black},
\begin{align}\label{perimeter_c1}
1\le   C(A+1)(1+L^3)\left( \rVert \partial_1\Delta^{-1}(1_D-1_{D^*})- \Omega\rVert_{\dot{H}^1}+\rVert \Omega\rVert_{L^2}\right),
\end{align}
where  $D^*:=\left\{ (x_1,x_2)\subset \mathbb{R}\times \mathbb{R}^-: (x_1,-x_2)\in D\right\}$.
\end{lemma}
\begin{remark}
As explained in Remark~\ref{remark_zeromean}, if we simply choose $\Omega=\partial_1\Delta^{-1}(1_D-1_{D^*})$ in \eqref{perimeter_c1}, then we get
\[
1\le C(A+1)(1+L^3)\rVert \partial_1\Delta^{-1}(1_D-1_{D^*})\rVert_{L^{2}}.
\]
This inequality indeed tells us that assuming $\int_{\mathbb{R}^2}1_{D}(x)x_2dx$ is bounded, the perimeter of $\partial D$ grows to infinity, as $\rVert \partial_1(1_D-1_{D^*})\rVert_{\dot{H}^{-2}}$ goes to zero.  In our proof of the main theorem, we will use the slightly finer estimate stated in the lemma, due to the lack of smallness of $\rVert \partial_1\rho(t)\rVert_{\dot{H}^{-2}}$.
\end{remark}
\begin{proof}
Without loss of generality, let us assume that $\inf\left\{ x_1: (x_1,x_2)\in D\right\} =0$ so that  a far-right point of $\partial D$ can be denoted by $(L,x_2^r)$ for some $x_2^r>0$. Using that $|D|=1$, we see 
\begin{align}
\int_{\mathbb{R}^2,\ \left\{x_{2}>4A\right\}}1_Ddx &\le \frac{1}{4A} \int_{\mathbb{R}^2,\ \left\{x_{2}>4A\right\}}1_D(x)x_2dx\le \frac{1}{4},\label{est12}\\
\int_{\mathbb{R}^2,\ \left\{x_{2}<\frac{1}{4L}\right\}}1_Ddx &=\int_{-\infty}^{\infty}\int_{0}^{\frac{1}{4L}}1_D(x)dx\le \frac{1}{4}\label{est13}.
\end{align}
This implies  
\begin{align}\label{AandL}
4A> \frac{1}{4L}.
\end{align}
Indeed, if it were not true, we would have 
\[
1 = |D| = \int_{\mathbb{R}^2}1_Ddx \le \int_{\mathbb{R}^2,\ \left\{x_{2}>4A\right\}} 1_D dx + \int_{\mathbb{R}^2,\ \left\{x_{2}<\frac{1}{4L}\right\}}1_Ddx \le \frac{1}{2},
\]
which is a contradiction. Moreover, the above estimates \eqref{est12} and \eqref{est13} imply that  at least a half of $D$ is contained in the horizontal strip bounded by $\left\{x_2= 4A\right\}$ and $\left\{ x_2=\frac{1}{4L}\right\}$. Therefore,  taking away $\left\{0\le x_1\le \frac{1}{32A},\ 0\le x_2\le 4A\right\}\cup\left\{L-\frac{1}{32A}\le x_1\le L,\  0\le x_2\le 4A\right\}$ from $D$, whose total measure is at most $\frac{1}{4}$,  we have that  
\begin{align}\label{confined_D}
\left|\left\{ x\in D: \frac{1}{32A} \le x_1 \le L-\frac{1}{32A},\quad \frac{1}{4L}\le x_2\le4A\right\} \right| \ge \frac{1}{4}.
\end{align}

Towards the proof of the lemma, we choose nonnegative smooth functions $g(x_1)$ and $h(x_2)$ satisfying
\begin{align*}
\begin{cases}
\text{$g(x_1)=0$ if $x_1\le 0$},\\
\text{$0\le g'(x_1)\le 1$ for $x_1\in (0,L)$ and $g'(x_1)=1$ for $x_1\in (\frac{1}{32A},L-\frac{1}{32A})$},\\
\text{$g(x_1)$ is symmetric about the axis $\left\{x_1=L\right\}$, that is, $g(x_1)=g(2L-x_1)$ for $x_1\ge L$},\\
\text{$0\le g(x_1)\le L$, $|g'(x_1)|\le 1$ and $|g''(x_1)|\le 32A$ for all $x_1\in\mathbb{R}$,}
\end{cases}
\end{align*}
and
\begin{align*}
\begin{cases}
\text{$h(x_2)= 0$ for $x_2\le 0$ or $x_2\ge 4A+\frac{1}{4L}$},\\
\text{$h(x_2) = 1$ for $x_2\in (\frac{1}{4L},4A)$},\\
\text{$0\le h(x_2)\le 1$, $|h'(x_2)|\le 4L$ and $|h''(x_2)|\le 32L^2$ for all $x_2\in\mathbb{R}$.}
\end{cases}
\end{align*}
A construction of $g,h$ satisfying above properties is straightforward. For such $g,h$, we define 
\begin{align*}
f(x_1,x_2):=g(x_1)h(x_2) \text{ for $x_2\ge 0$ and } f(x_1,x_2)=-f(x_1,-x_2), \text{ for $x_2< 0$.}
\end{align*}
Clearly, $|\text{supp}(f)|\le C(AL+1)$, thus the above properties of $g,h$ give us that 
\begin{align}\label{support_f2}
|\text{supp}(\nabla f)|\le C(AL+1)\le CAL,
\end{align}  
where the last inequality is due to \eqref{AandL}.
Furthermore, noticing that $g(x_1)$ is linear for $x_1\in (\frac{1}{32A},L-\frac{1}{32A})$, it is not difficult to see that
\begin{align}\label{support_f3}
|\text{supp}(\Delta f)|\le C.
\end{align}

Next, denoting $\mu:=1_D - 1_{D^*}$ and following the same computations in \eqref{l2dual}, we get  
 \begin{align}\label{fract}
{\int_{\mathbb{R}^2} \mu(x)\partial_1f(x)dx}\le \rVert \partial_1\Delta^{-1}\mu - \Omega\rVert_{\dot{H}^1}\rVert \nabla f\rVert_{L^2} +\rVert \Omega\rVert_{L^2}\rVert \Delta f\rVert_{L^2},
 \end{align}
for any $ \Omega$ satisfying \eqref{omega_relaxation}\color{black}.
Using $x_2$-odd symmetry of $\mu$ and $f$, we see that $\int_{\mathbb{R}^2}\mu(x)\partial_1 f(x)dx = 2\int_{D}g'(x_1)h(x_2)dx$, while the properties of $g,h$ give us that
\begin{align*}
\int_{D}g'(x_1)h(x_2)dx\ge \int_{\frac{1}{32A}}^{L-\frac{1}{32A}}\int_{\frac{1}{4L}}^{4A}1_D(x_1,x_2)dx_2dx_1\ge \frac{1}{4},
\end{align*}
where the last inequality follows from \eqref{confined_D}. Thus, we have
\begin{align}\label{nue}
\int_{\mathbb{R}^2}\mu(x)\partial_1 f(x)dx \ge \frac{1}{2}.
\end{align} On the other hand, we have that
\[
\rVert \Delta f \rVert_{L^\infty}\le \rVert g''\rVert_{L^\infty} + \rVert h''\rVert_{L^\infty}\rVert g\rVert_{L^\infty} \le 32A + 32L^3,\quad \rVert \nabla f\rVert_{L^\infty}\le \rVert g'\rVert_{L^\infty}+\rVert h'\rVert_{L^\infty}\rVert g\rVert_{L^\infty}\le 1+4L^2.
\]
Combining this with \eqref{support_f2} and \eqref{support_f3}, we get
\[
\rVert \Delta f\rVert_{L^2}\le C(A+L^3)\le C(A+1)(1+L^3),\quad \rVert \nabla f \rVert_{L^2} \le C(AL^3 +AL )\le C(A+1)(1+L^3).
\]

Plugging this and \eqref{nue} into \eqref{fract}, we obtain the desired estimate \eqref{perimeter_c1}.
\end{proof}

\section{Proof of the main theorem}
In this section, we prove the main theorem of the paper. Let $\rho(t)=1_{D_t}-1_{D_t^*}$ be the global solution with the initial data $(\rho_0,u_0)$ satisfying the assumptions \ref{assumption1} and \ref{assumption2}. In the rest of the proof, $C$ denotes some positive constant that depends on only $(\rho_0,u_0,\nu)$ and might vary from line to line. 

  From \eqref{bound2} and Lemma~\ref{uniform_vs}, we have that for any $n\in\mathbb{N}$, we can find $T_n>0$ such that $T_n\mapsto \infty$ and 
\[
\int_{T_n}^{2T_n} \rVert \omega(t)\rVert_{L^2}^2 + \rVert \partial_1 \Delta^{-1}\rho(t) - \nu\omega(t)\rVert_{\dot{H}^1}^2 dt \le \frac{1}{n}.
\]
Therefore, there exists $t_n\in [T_n,2T_n]$ such that 
\begin{align}\label{t_nomega}
\rVert \omega(t_n)\rVert_{L^2} + \rVert \partial_1 \Delta^{-1}\rho(t_n) - \nu\omega(t_n)\rVert_{\dot{H}^1} \le \frac{1}{\sqrt{n T_n}}.
\end{align}
We prove the growth of curvature \ref{theorem1} first.
\begin{proof}[\textup{\textbf{Proof of Theorem~\ref{curvature}, part (a)}}] 
 Let $B_r(x^*)$ be the largest disk contained in $D_{t_n}$, centered at $x^* = (x_1^*, x_2^*)$ with radius $r$. Then the Pestov--Ionin theorem (see \eqref{pIlem}) tells us that
\begin{align}\label{radicurv}
r\ge \frac{1}{\max_{x\in\partial D_{t_n}}|\kappa(t_n)|},
\end{align}
where $\kappa(t_n)$ is the signed curvature of $\partial D_{t_n}$.  Applying Lemma~\ref{curvature_lem} with $D=D_{t_n}$ and $\Omega=\omega(t_n)$, we obtain
\[
r^3\le {C}\left( \rVert \partial_1\Delta^{-1}\rho(t_n) - \omega(t_n)\rVert_{\dot{H}^1} +\rVert \omega(t_n)\rVert_{L^2}\right)\le \frac{C}{\sqrt{nT_n}},
\]
where the last inequality follows from \eqref{t_nomega}.  Combining this with \eqref{radicurv} yields that
\[
\max_{x\in\partial D_{t_n}}|\kappa(t_n)| \ge (nT_n)^{\frac{1}{6}} \ge (nt_n)^{\frac{1}{6}},
\]
where we used $t_n\in [T_n,2T_n]$ for the last inequality. Since $t_n\mapsto \infty$, as $n\to \infty$, we obtain the desired infinite growth of the curvature.
\end{proof}

\begin{proof}[\textup{\textbf{Proof of Theorem~\ref{curvature}, part (b)}}]
Let $L_n>0$ be the distance between a far-left and a far-right points on $\partial D_{t_n}$.  From the energy conservation \eqref{energy_conservation}, we see that $E_P(t)\le E_T(0)<C$ for all $t>0$. Therefore,  applying Lemma~\ref{periLem} with $D=D_{t_n}$, $\Omega=\omega(t_n)$ and $A=\frac{1}{2}E_P(t_n)\le C$, we obtain
\[
1\le C(L_n^3+1)\left( \rVert \partial_1\Delta^{-1}\rho(t_n)-\omega(t_n)\rVert_{\dot{H}^1} + \rVert \omega(t_n)\rVert_{L^2}\right)\le \frac{C}{\sqrt{nT_n}}(L_n^3+1),
\]
where the last inequality follows from \eqref{t_nomega}.
Since $T_n\to \infty$ as $n\to \infty$, the above inequalities give us that $L_n \ge C (nT_n)^{\frac{1}{6}}\ge C(nt_n)^{\frac{1}{6}}$ for all $n\in\mathbb{N}$, where $t_n\le 2T_n$ was used to justify the second inequality. Since $D_{t_n}$ is simply connected, this proves the desired infinite growth of perimeter.
\end{proof}

\subsection*{Acknowledgements}
The author was partially supported by the SNF grant 212573-FLUTURA and the Ambizione fellowship project PZ00P2-216083. The author also extends gratitude to  Eduardo Garc\'ia-Ju\'arez and Yao Yao for their valuable discussions and insightful suggestions during the course of this research.
\bibliographystyle{abbrv}
\bibliography{references}

\end{document}